\documentclass[a4paper,12pt]{amsart}

\usepackage{amssymb,mathrsfs,bm}
\usepackage[dvipsnames]{xcolor}
\usepackage{mathtools,galois,braket,esint}
\usepackage[utf8]{inputenc}
\usepackage[margin=0.5in]{geometry}
\usepackage{hyperref,hypcap}
\hypersetup{colorlinks=true,allcolors=black,bookmarksdepth=3}
\usepackage[inline]{enumitem}
\usepackage{array}
\usepackage[utf8]{inputenc}
\usepackage{tikz, tikz-3dplot}
\usepackage{amsmath, amsthm, amssymb,graphics }
\usetikzlibrary{arrows,shapes,positioning}
\usetikzlibrary{intersections}
\usetikzlibrary{decorations.pathreplacing,calligraphy}
\newcommand{\dd}{\mathop{}\!\mathrm{d}}

\usepackage{amsmath,amsthm,amssymb,amscd,color, xcolor,mathtools,url,tikz}
\usepackage{bbm}
\usepackage{pgfplots}

\theoremstyle{plain}
\newtheorem{theorem}{Theorem}[section]

\theoremstyle{definition}

\theoremstyle{remark}

\numberwithin{equation}{section}
\numberwithin{figure}{section}
\numberwithin{table}{section}
\allowdisplaybreaks

\begin{document}

\title[A review on asymptotic stability of solitary waves]{A review on asymptotic stability of solitary waves in nonlinear dispersive problems in dimension one}

\begin{abstract} We review asymptotic stability of solitary waves for nonlinear dispersive equations set on the line. Our focus is threefold: first, the nonlinear Schr\"odinger equation; second, the notion of full asymptotic stability (which states that perturbations of a solitary wave decompose globally into a solitary wave and a decaying solution); and third, spectral methods. Besides this focus, we summarize the state of the art in a broader context, including nonlinear Klein-Gordon equations, the notion of local asymptotic stability, and virial methods.
\end{abstract}

\author[P. Germain]{Pierre Germain}
\address{Pierre Germain, Department of Mathematics, Huxley Building, South Kensington Campus,
Imperial College London, London SW7 2AZ, United Kingdom}
\email{pgermain@ic.ac.uk}

\maketitle

\setcounter{tocdepth}{1}

\tableofcontents

This text aims at reviewing the asymptotic stability of solitary waves in  nonlinear dispersive equations set in Euclidean space. Given the breadth of the subject, we had to emphasize some directions rather than others, and our choice was the following
\begin{itemize}
\item We set the \textit{space dimension equal to one} - this is the best understood case as of now.
\item We focused on the \textit{nonlinear Schr\"odinger equation} \eqref{NLS}, which is one of the canonical models.
\item We favored \textit{full asymptotic stability} over various other notions of stability which will be reviewed below. It can be roughly defined as follows: perturbations of a solitary wave decompose globally, as time goes to infinity, into the sum of a solitary wave and a pointwise decaying part.
\item  Finally, \textit{spectral methods} seem unavoidable to prove that this strong notion of stability is satisfied; here, spectral is understood in the sense of the Fourier transform, or of the linearized group around a traveling wave.
\end{itemize}

We believe that the ideas and methods which we will present for \eqref{NLS} are valid for most nonlinear dispersive problems, in dimension one and higher - indeed, we will see that they have immediate counterparts for nonlinear Klein-Gordon equations.
Besides full asymptotic stability, local asymptotic stability has seen important progress over the last decade, in connection with viral or energy methods. We will describe the results and sketch the associated ideas.

Finally, we refer to the textbooks \cite{Pava,SulemSulem} and the reviews \cite{CuMa,KowalczykMartelMunoz1bis,MartelICM,Martel2024,Soffer,Tao2009} to give the reader a more complete picture of the subject.

\section{The equation, its solitary waves, and their stability}

\subsection{The equation and its solitary waves}
We consider the Cauchy problem for the nonlinear Schr\"odinger equation
\begin{equation} \tag{NLS}
\label{NLS} 
\begin{cases} 
i \partial_t v - \partial_x^2 v - F'(|v|^2) v= 0 \\
v(t=0) = v_0.
\end{cases}
\end{equation}
set on the real line $v=v(t,x)$, $t,x \in \mathbb{R}$.
It is stemming from the Hamiltonian
$$
H(v) = \int [|\partial_x v|^2 - F(|v|^2)] \dd x.
$$
The interaction potential $F$ will be assumed to be smooth and to have a non-degenerate local minimum at zero\footnote{Besides smoothness, all we need of $F$ is the existence of stationary waves, which is a consequence of having a local minimum at zero, but not equivalent, see \cite{BL} for the exact condition.}. Stationary waves of the type
$$
v(t) = e^{-it\omega} \Phi_\omega, \qquad \omega>0,
$$
are given by solutions of 
\begin{equation} \label{eq:soliton}
\partial_x^2 \Phi_\omega - \omega \Phi_\omega + F'(\Phi_\omega^2) \Phi_\omega = 0.
\end{equation}
Under our assumptions on $F$, there exists a unique solution up to translation of this ODE on a non-trivial interval $\omega \in (0,\omega^*)$.
	
For $p, \gamma, y \in \mathbb{R}$, Galilean, phase and translation symmetries
$$
v(t,x) \mapsto e^{i(p x + p^2 t + \gamma)} v(t,x+2p  t - y)
$$
leave the set of solutions of~\eqref{NLS} invariant. In particular, this gives the family of traveling waves
\begin{equation} \label{travelingfamily}
e^{i( p x +( p^2-\omega) t + \gamma)} \Phi_\omega (x+2p t - y).
\end{equation}

\subsection{Existence and uniqueness of solutions} These questions will not concern us in the present review paper. Suffice it to say that local well-posedness holds in $H^1$ and that the solutions are global if $F$ grows slower than at a cubic rate at infinity ($|F(z)| \lesssim |z|^{3-\delta}$ with $\delta>0$). If $F$ grows faster than cubic, then blowup becomes possible, even for perturbations of solitary waves. While this has to be kept in mind, the focus of the present paper will be the opposite case, where solitons are stable.

For this and more, we refer to the classical textbooks \cite{Cazenave,SulemSulem,Tao}.

\subsection{Different kinds of stability} Whether the solitary waves $\Phi_\omega$ are stable or not is of foremost importance from a mathematical as well as a physical viewpoint, since stable objects are of greater relevance to the dynamics. Stability will depend on the norm $\| \cdot \|$ under consideration - it will not specified for the time being. Once the topology is chosen, different kinds of stability can be considered

\medskip

\textit{Lyapunov stability} would be asking for any $\epsilon >0$ the existence of $\delta >0$ such that $\sup_{t>0} \| v(t) - e^{it} \Phi_\omega \| < \epsilon$ provided $\| v_0 - \Phi_\omega \| < \delta$. As is well known, this is too naive and does not hold. Indeed, the data $\Phi_\omega$ gives the solution $e^{-i \omega t} \Phi_\omega(x)$; modifying this data slighlty to $\Phi_\omega e^{ipx}$ (for any $p \neq 0$) leads to the solution $e^{i(px + p^2 t)} \Phi_\omega(x +2pt)$ which inexorably drifts away from the earlier solution.

\medskip

\textit{Orbital stability} is meant to fix the problem that was just identified by incorporating the symmetries of the equation. Given $\epsilon>0$, we now look for $\delta>0$ such that, if $\| v_0 - \Phi_\omega \| < \delta$, then
\begin{equation}
\label{orbitalstability}
\sup_{t >0} \, \inf_{\omega,p,\theta,y \in \mathbb{R}} \, \| v(t,x) - e^{i(px + \theta)} \Phi_\omega(x + y) \| < \epsilon.
\end{equation}
Orbital stability is by now completely understood \cite{CL,GSS2,Weinstein};  it holds in the $H^1$ topology if and only if
\begin{equation}
\label{defcomega}
c_\omega = \frac{d}{d\omega} \int | \Phi_\omega |^2 \dd x >0
\end{equation}
(ignoring the degenerate case where this quantity is zero). This is sometimes called the Vakhitov-Kolokolov condition.

\medskip

\textit{Local asymptotic stability} describes a stronger property than orbital stability: namely, it is asking that the solutions converge locally to the soliton. For a cutoff function $\chi$, we ask that, if $\| v_0 -  \Phi_\omega \|$ is sufficiently small, then
$$
\sup_{R>0} \, \lim_{t \to \infty} \, \inf_{\omega,p,\theta,y \in \mathbb{R}} \left\| \chi \left( \frac {x+y} R \right) \left[ v(t,x) - e^{i(px + \theta)} \Phi_\omega(x + y) \right] \right\| = 0.
$$

\medskip

\textit{Full asymptotic stability} is asking for a full description of the solution if $\| v_0 - \Phi_\omega \|$ is sufficiently small. Namely, we are asking for a decomposition
\begin{equation} \label{decfull}
v(t,x) = e^{i( p(t)x + \theta(t))} \Phi_{\omega(t)} (x + y(t)) + \{\mbox{time-decaying solution} \} \qquad \mbox{as $t \to \infty$}
\end{equation}
where the asymptotic behavior of $p(t),\theta(t),\omega(t),y(t)$ is described and the decaying solution, or radiation, $w(t)$ will be characterized as undergoing scattering or modified scattering. 

\medskip

\textit{The soliton resolution conjecture} is asking for a decomposition similar to \eqref{decfull} with two important differences: first, it should hold for all (or almost all, in an appropriate sense) data, and second, it shoud allow for a finite number of solitary waves on the right-hand side.
Such a statement seems out of reach of present tools for \eqref{NLS}, except in the completely integrable case that we will come back to; but full asymptotic stability is the first step towards this much more ambitious goal.

\medskip

Having recapitulated these different notions of stability, we now formulate the question that will be at the heart of the present review.

$$
\boxed{\mbox{Question: For which $F$ can full asymptotic stability be established?}}
$$

\subsection{Notation}
We adopt the following normalization for the Fourier transform of $u$
$$
\widehat{u}(\xi) = \frac{1}{\sqrt{2\pi}} \int_{-\infty}^\infty u(x) e^{- i x \xi} \dd \xi, \qquad u(x) = \frac{1}{\sqrt{2\pi}} \int_{-\infty}^\infty \widehat{u}(\xi)  e^{ i x \xi} \dd \xi.
$$
Sobolev and weighted $L^2$ spaces are denoted as follows
$$
\| u \|_{H^s} = \| \langle \partial_x \rangle^s u \|_{L^2} \qquad \| u \|_{L^{{2,s}}} = \| \langle x \rangle^s u \|_{L^2}.
$$

\section{The main characters}

In this section, we want to present in greater detail the two main characters of the story: on the one hand, the solitary wave, and on the other, the decaying wave or radiation, whose behavior is characterized by scattering or modified scattering. As we saw earlier, full asymptotic stability is asking for a behavior as $t\to \infty$ which is a linear combintion of these two types of solutions.

We will also mention a third character, namely breathers. These objects do not resolve as $t \to \infty$ into a sum of solitary wave and radiation, and thus we should be careful to avoid them given our definition of asymptotic stability.

\subsection{The solitary waves} 
\label{subsectionsolitary}
Recall that stationary waves
$$
v(t) = e^{-it\omega} \Phi_\omega 
$$
are given by solutions of the ODE
\begin{equation} 
\partial_x^2 \Phi_\omega - \omega \Phi_\omega + F'(\Phi_\omega^2) \Phi_\omega = 0.
\end{equation}
By our assumptions on $F$, there is a unique solution up to translation on an interval $\omega \in (0,\omega^*)$, for some $\omega^*>0$. Furthermore, $\Phi_\omega$ is even, positive, decreasing on $x>0$, and exponentially decreasing at infinity, along with its derivatives (see~\cite{BL} for these facts and optimal conditions on $F$ for the existence of stationary waves).

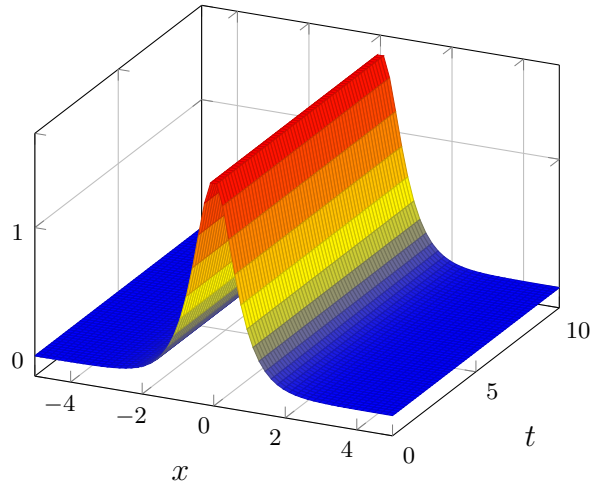
\begin{figure} \label{pictureSW}
\begin{tikzpicture}
\begin{axis} [
    xlabel = $x$, ylabel = $t$,
    ticklabel style = {font = \scriptsize},
	grid
]
\addplot3 [surf, domain=-5:5, y domain=0:10, samples=60] 
	{ 1.6/cosh(2*x) };
\end{axis}
\end{tikzpicture}
\caption{The modulus of the stationary wave solution of the cubic \eqref{NLS}: $\Phi_1(x) = |e^{it} \Phi_1(x)|$}
\end{figure}

In general, there does not exist an explicit formula for the solitary waves, but this is for instance the case if $F$ is a pure power: 
$$
\mbox{if $F(z) = \frac{1}{\sigma +1} z^{\sigma +1}$}, \qquad \Phi_1 (x) = \frac{(\sigma +1)^{\frac{1}{2\sigma}}}{\cosh(\sigma x)^{\frac{1}{\sigma}}}.
$$
The solitary wave is shown in Figure \ref{pictureSW}.

Still in the pure power case, the equation enjoys a scaling invariance which leads to the formula
$$
\Phi_\omega(x) = \omega^{\frac{1}{2\sigma}} \Phi_1 (\omega^{\frac 12} x).
$$

This formula has two interesting consequences. First, substituting it in the criterion \eqref{orbitalstability} for orbital stability gives the condition $\sigma < 2$ for orbital stability. Second, if $\sigma= 1$ (cubic \eqref{NLS}), the above formula yields
$$
\| \Phi_\omega \|_{H^1} \sim \omega^{\frac 14} \qquad \mbox{and} \qquad \| \Phi_\omega \|_{L^{2,1}} \sim \omega^{-\frac 14}
$$
(the same would be true if $F(z) $ is quadratic to leading order as $z \to 0$, rather than exactly equal to $z^2$). This shows the importance of the norm in which the perturbation to the solitary wave is measured:
\begin{itemize}
\item Any $H^1$ neighborhood of a given solitary wave will contain small solitary waves, which invalidates full asymptotic stability as stated in \eqref{decfull} (of course, solitary waves can only be added rigorously in the asymptotic regime $t\to \infty$, see \cite{MartelMerle} for a construction of a multi-soliton solution). The obvious fix is to add small solitons in \eqref{decfull}, but proving such a decomposition might be as hard as the soliton resolution conjecture.
\item But a neighborhood in $L^{2,1}$ of a given solitary wave should exclude the presence of small solitons in the background.
\end{itemize}

\subsection{Time-decaying solutions} We now turn to small data $u_0$ (say in $H^1 \cap L^{2,1}$) which lead to decaying solutions, thus excluding the appearance of solitary waves. 

\subsubsection*{Scattering} 
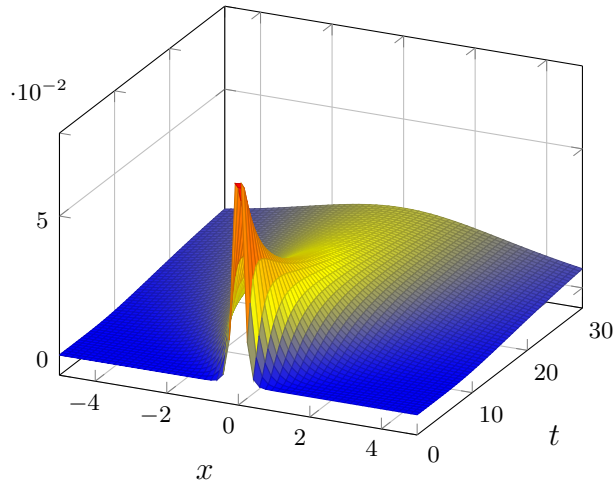
\begin{figure} \label{picturescattering}
\begin{tikzpicture}
\begin{axis} [
    xlabel = $x$, ylabel = $t$,
    ticklabel style = {font = \scriptsize},
	grid
]
\addplot3 [surf, domain=-5:5, y domain=0:30, samples=60] 
	{ .1 * (y+1)^(-.5) / cosh(10*x/(1+y)) };
\end{axis}
\end{tikzpicture}
\caption{The modulus of a small scattering solution}
\end{figure}

If the potential $F$ vanishes to order $>2$ at the origin \footnote{In other words, $F(z) = O(z^3)$, since we are assuming $F$ to be smooth, but more generally, $F(z) = O(z^\alpha$ with $\alpha>2$ would suffice.} the solutions scatter, which means that they behave as $t \to \infty$ like linear solutions, or in other words that the  \textit{profile} $f$ converges:
$$
f(t) = e^{it \partial_x^2} v(t) \rightarrow f_\infty \qquad \mbox{as $t \to \infty$.}
$$
In terms of the solution $v$, it resembles asymptotically $e^{-it\partial_x^2} f_{\infty}$, which, by the stationary phase lemma, implies that
$$
u(t,x) \asymp e^{it\partial_x^2} f_{\infty} (x) \asymp e^{- i \frac \pi 4} {\frac 1 {\sqrt{2t}}}  \widehat{f_\infty} \left( \frac { x} {2t} \right) e^{i \frac{x^2}{4t}} \qquad \mbox{as $t \to \infty$}.
$$
A scattering solution is represented in Figure \ref{picturescattering}.

\subsubsection*{Modified scattering} If the potential $F$ vanishes to order $2$ exactly at the origin (which is the "generic" situation if $F$ has a minimum at zero), then the nonlinearity is cubic to leading order. This is critical from the point of view of asymptotic behavior: namely, if $u$ satisfies the standard decay estimates for linear solutions $\| u \|_{L^2} \sim 1$, $\| u (t) \|_{L^\infty} \sim t^{-\frac 12}$, a naive computation gives logarithmic growth
\begin{equation}
\label{loggrowth}
\int_0^T \| |u(t)|^2 u(t) \|_{L^2} \dd t \lesssim \int_0^T \| |u(t) \|_{L^\infty}^2 \| u(t) \|_{L^2} \dd t \sim \int_0^T \frac{\dd t}{t} \sim \log T \qquad \mbox{as $T \to \infty$}.
\end{equation}
This will ultimately lead to a correction to scattering.

To understand the nature of this correction, we will follow the general idea of the space-time resonance method and view the nonlinear term as an oscillatory integral in Fourier space. The equation \eqref{NLS} can be written as as an evolution problem for the profile $f(t)$ in Fourier space: using basic properties of the Fourier transform, we obtain
\begin{equation}
\label{equationherisson}
i \partial_t \widehat{f}(t,\xi) = \frac{c}{\pi} \iint e^{it (\xi^2 - \eta^2 - \sigma^2 - (\xi-\eta-\sigma)^2)} \widehat{f}(t,\eta) \widehat{f}(t,\sigma) \overline{\widehat{f}(t,-\xi+\eta+\sigma)} \dd \eta \dd \sigma + \{ \mbox{negligible} \}
\end{equation}
(here, $c$ is the coefficient of $z^2$ in the expansion of $F(z)$ at zero, or equivalenty the coefficient of $|u|^2 u$ in \eqref{NLS}; as already discussed, higher order terms do not matter here).
To deal with this integral, we assume that $\widehat{f}$ is smooth in $t$ and $\xi$, making the application of the stationary phase lemma legitimate. The phase is
$$
\varphi(\eta,\sigma) = \xi^2 - \eta^2 - \sigma^2 + (\xi-\eta-\sigma)^2 = 2( \xi^2 - \xi \eta - \xi \sigma + \eta \sigma).
$$
It is stationary if 
$$
\partial_\eta \varphi(\xi,\eta,\sigma) = \partial_\sigma \varphi(\xi,\eta,\sigma) = 0 \qquad \Longleftrightarrow \qquad \eta = \sigma = \xi.
$$
At the stationary point, $\phi(\xi,\xi,\xi) = 0$ and the Hessian is $\begin{pmatrix} 0 & 2 \\ 2 & 0 \end{pmatrix}$ with determinant $4$ and signature zero. By the stationary phase lemma, the right-hand side of \eqref{equationherisson} is $\asymp \frac{c}{t} ||\widehat{f}(t,\xi)|^2 \widehat{f}(t,\xi)$, leading to the following asymptotic ODE for $\widehat{f}$
\begin{align*}
i \partial_t \widehat{f}(t,\xi) = \frac{c}{t} |\widehat{f}(t,\xi)|^2 \widehat{f}(t,\xi).
\end{align*}
Integrating this ODE gives a logarithmic correction to scattering, namely 
$$
\widehat{f}(t,\xi) \asymp \widehat{f_\infty}(\xi)e^{ic |\widehat{f_\infty}(\xi)|^2 \log t } \qquad \mbox{or} \qquad 
u(t,x) \asymp e^{- i \frac \pi 4} {\frac 1 {\sqrt{2t}}}  \widehat{f_\infty} \left( \frac { x} {2t} \right) e^{i \left[ \frac{x^2}{4t} + c \left| \widehat{f_\infty} \left(\frac { x} {2t} \right) \right|^2 \log t \right]} \qquad \mbox{as $t \to \infty$}.
$$
for an asymptotic profile $f_\infty$.

We followed the argument in Kato-Pusateri \cite{KP} in deriving heuristically this formula for modified scattering. We refer to that pater for further details and a full justification. Many other approaches were put forward to show modified scattering \cite{HN,IT,LS}, see the review \cite{Murphy}.

\subsection{Breather solutions} \label{subsectionbreathers}
\begin{figure} \label{pictureSY}
\begin{tikzpicture}
\begin{axis} [
    xlabel = $x$, ylabel = $t$,
    ticklabel style = {font = \scriptsize},
	grid
]
\addplot3 [surf, domain=-6:6, y domain=0:160, samples=60] 
	{ 4*1.6*((cosh(3*x) + 3* cos(8*y) * cosh(x))^2 + (3* sin(8*y) * cosh(x))^2)^(.5)/(cosh(4*x) + 4 * cosh(2*x) + 3*cos(8*y)) };
\end{axis}
\end{tikzpicture}
\caption{The modulus of the Satsuma-Yajima breather}
\end{figure}
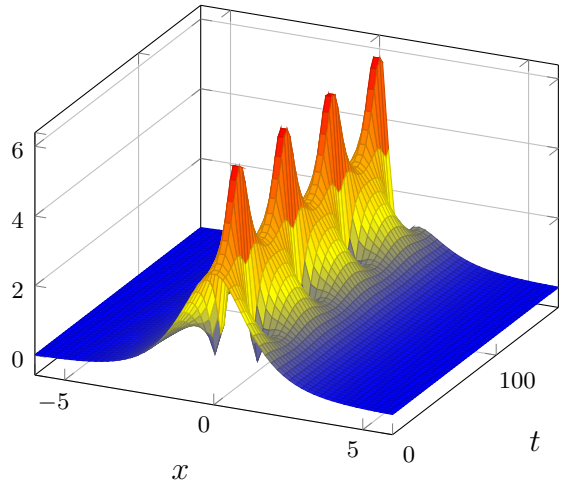

We borrow the following formula from \cite{AlejoFanelliMunoz}, see also \cite{AA,Olmedilla}. For any $c_1,c_2 >0$, if $\gamma_{\pm} = c_1 \pm c_2$, then an explicit solution of the cubic \eqref{NLS} (ie $F(z) = z^2$) is 
$$
B(t,x) = \frac{2\sqrt 2 \gamma_+ \gamma_- e^{-i c_1^2 t}[ c_1 \cosh (c_2 x) + c_2 e^{i \gamma_+ \gamma_- t} \cosh(c_1 x)]}{\gamma_-^2 \cosh(\gamma_+ x) + \gamma_+^2 \cosh(\gamma_- x) + 4 c_1 c_2 \cos(\gamma_+ \gamma_- t)}.
$$
It is instructive to consider a few special cases for $c_1$ and $c_2$
\begin{itemize}
\item If $c_1=1$, $c_2=0$, this is the solitary wave $e^{-it} \Phi_1$.
\item If $c_1=1$, $c_2=3$, this is the Satsuma-Yajima breather \cite{SY} which is depicted in Figure \ref{pictureSY}
$$
\frac{4 \sqrt 2 e^{-it} (\cosh(3x) + 3 e^{8it} \cosh x)}{\cosh(4x) + 4 \cosh(2x) + 3 \cos(8t)}.
$$
\item If $c_1 =1$, $c_2 \ll 1$, one obtains a breather which looks like the solitary wave $e^{-it} \Phi_1$ for $x$ small, but which decays much slower, like $c_2 e^{-c_2|x|}$, as $x \to \infty$.
\end{itemize}

As a consequence of the last point, any neighborhood of $\Phi_1$ in $H^1$ contains nonlinear solutions which are time-periodic, and are not  resolved into solitary wave plus radiation as time goes to infinity - even though these breather solutions are unstable. This is clearly an obstacle to full and even local asymptotic stability as stated above. This obstacle can be lifted in two ways: either by considering a stronger topology than $H^1$, or by replacing the family of solitons in the definition of asymptotic stability by the larger family including breathers - but a practical implementation of this idea seems very delicate!

\section{Linearization around a solitary wave and decay}

In order to linearize atound the stationary wave $e^{-it\omega} \Phi_\omega$, one filters first the time oscillation by setting $u(t) = e^{it\omega} v(t)$. The equation satisfied by $u$ is now
$$
i \partial_t u + \omega u - \partial_x^2 u - F'(|u|^2) u = 0.
$$
For this new equation, $\Phi_\omega$ is a stationary solution. Linearizing around it leads to a problem which is not complex-linear in $u$ (because of the nonlinearity). To make this equation complex-linear, we write it as a vector problem in $U = \begin{pmatrix} u \\ \overline{u} \end{pmatrix}$. The resulting linearized problem is
$$
i \partial_t U + \mathcal{H}_\omega U = 0.
$$
where 
\begin{align*}
&   \mathcal{H}_\omega =  \begin{pmatrix}
-\partial_x^2 + \omega & 0  \\
0 &  \partial_x^2 - \omega
\end{pmatrix} - \begin{pmatrix}
V_+  & V_-
\\
-V_- & -V_+
\end{pmatrix} \\
& V_- = F''(\Phi^2) \Phi^2 \qquad \mbox{and} \qquad V_+=F'(\Phi^2)+V_-.
\end{align*}
The operator $\mathcal{H}_\omega$ is a vector Schr\"odinger operator, with a matrix potential; such operator shares many common features with the more classical scalar Schr\"odinger operators, but there are also important differences, first of all since it is not self-adjoint.

\subsection{Spectrum of the linearized operator}
The spectrum of $\mathcal{H}_\omega$ is obviously a fundamental ingredient to understand the stability of solitary waves. It can be understood rather precisely \cite{BP,ErdoganSchlag}
\begin{itemize}
\item \textit{The essential spectrum} equals $(-\infty,-\omega] \cup [\omega,\infty)$ by the Weyl criterion.

\medskip

\item \textit{Eigenvalues} belong to $\mathbb{R} \cup i \mathbb{R}$. 

\medskip

\item \textit{Unstable modes} are by definition eigenvalues in $i \mathbb{R} \setminus \{0\}$; they are absent if $c_\omega>0$.

\medskip

\item \textit{The generalized kernel} has geometric dimension at least $2$ and algebraic dimension at least $4$, which is the dimension of the manifold of solitary waves \eqref{travelingfamily}. Indeed, differentiating the equation with respect to the parameters $(p,\omega,\gamma,y)$ leads to the following four elements in the generalized kernel: $\begin{pmatrix} \Phi  \\ -\Phi \end{pmatrix}$, $\begin{pmatrix} \partial_\omega \Phi  \\ \partial_\omega \Phi \end{pmatrix}$, $\begin{pmatrix} \partial_x \Phi  \\ \partial_x \Phi \end{pmatrix}$, $\begin{pmatrix} x \Phi  \\ -x \Phi \end{pmatrix}$. If $c_\omega > 0$, the generalized kernel is spanned by them.
\item \textit{Embedded eigenvalues} (eigenvalues contained in the essential spectrum) are conjectured not to exist in \cite{ErdoganSchlag}, but it seems that a proof is only known in the pure power case \cite{KS,Perelman}.

\medskip

\item \textit{Internal modes} are by definition eigenvalues in $(-\omega_0,\omega_0)$; see \cite{Rialland} for a characterization of potentials $F$ which are cubic to leading order and for which internal modes are absent.

\medskip

\item \textit{Edge resonances} (in other words, bounded generalized eigenfunctions corresponding to the eigenvalues $\pm \omega$) exist for certain $F$, see below.
\end{itemize}

Figures \ref{spectrum-a} and \ref{spectrum-p} illustrate the spectrum of the linearized operator around the solitary wave in two emblematic cases: the pure power noninearity and the cubic-quintic nonlinearity. In both cases, it is interesting to observe how the spectrum evolves with the parameter.

\setlength{\unitlength}{1cm}

\begin{figure}
\label{spectrum-p}
\begin{picture}(20,4.5)
\thinlines 
\put(0,2){\vector(1,0){4}}  
\put(2,0){\vector(0,1){4}}  

\put(5,2){\vector(1,0){4}}  
\put(7,0){\vector(0,1){4}}  

\put(10,2){\vector(1,0){4}}  
\put(12,0){\vector(0,1){4}}  

\put(15,2){\vector(3,0){4}}  
\put(17,0){\vector(0,3){4}}  

\put(2,2){{\color{red}\circle*{0.25}}}
\put(7,2){{\color{red}\circle*{0.25}}}
\put(12,2){{\color{red}\circle*{0.4}}}
\put(17,2){{\color{red}\circle*{0.25}}}

\put(1,2){{\color{red}\circle*{0.25}}}
\put(3,2){{\color{red}\circle*{0.25}}}

\put(6.5,2){{\color{red}\circle*{0.25}}}
\put(7.5,2){{\color{red}\circle*{0.25}}}

\put(17,3){{\color{red}\circle*{0.25}}}
\put(17,1){{\color{red}\circle*{0.25}}}

\linethickness{1mm}
\put(0,2){\color{red}\line(1,0){1}}
\put(3,2){\color{red}\line(1,0){.8}}

\put(5,2){\color{red}\line(1,0){1}}
\put(8,2){\color{red}\line(1,0){.8}}

\put(10,2){\color{red}\line(1,0){1}}
\put(13,2){\color{red}\line(1,0){.8}}

\put(15,2){\color{red}\line(1,0){1}}
\put(18,2){\color{red}\line(1,0){.8}}

\put(.8,.3){\kern3pt $\sigma = 1$ (cubic NLS)}
\put(5.8,.3){\kern3pt\texttt{$1 < \sigma < 2$ }}
\put(10.8,.3){\kern3pt $\sigma = 2$ (quintic NLS)}
\put(15.8,.3){\kern3pt\texttt{$\sigma>2$ }}

\end{picture}

\caption{The spectrum of $\mathcal{H}_\omega$ for the pure power case $F(z) = z^{\sigma+1}$ after \cite{CGNT,ColesGustafson,KS}. Spectrum is in red, lines represent essential spectrum and points discrete spectrum or edge resonances. For $\sigma =1$, edge resonances are present; for $\sigma < 1 < 2$ internal modes; for $\sigma =2$ the generalized kernel degenerates, hence the bigger point at the origin; for $\sigma>2$ unstable modes appear. Observe how the eigenvalues bifurcate from the resonance for $\sigma \to 1$.}
\end{figure}
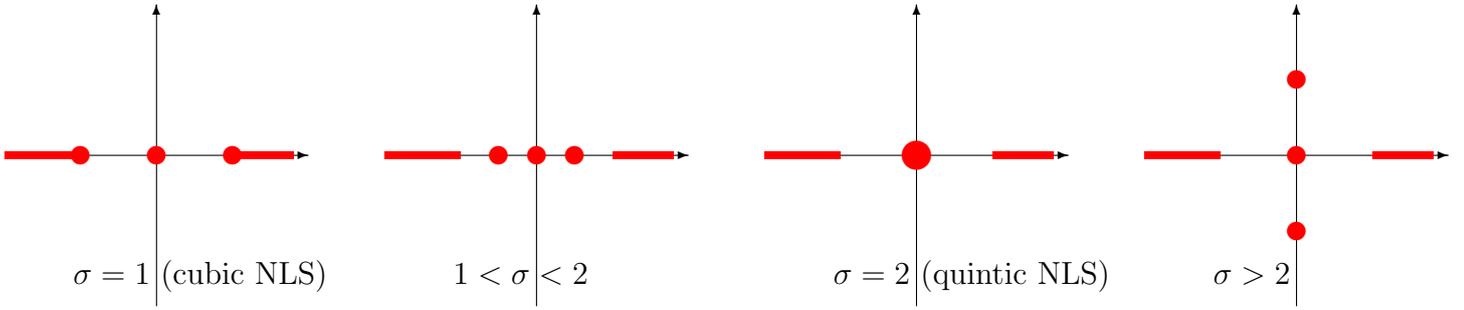

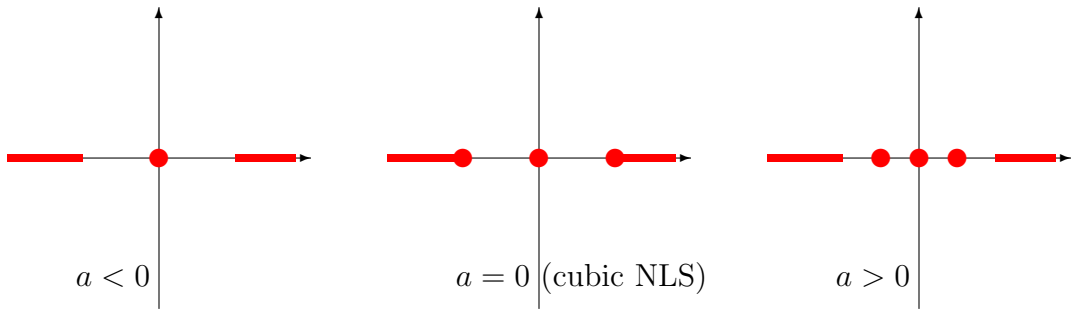
\begin{figure} \label{spectrum-a}
\begin{picture}(15,4.5)
\thinlines 
\put(0,2){\vector(1,0){4}}  
\put(2,0){\vector(0,1){4}}  

\put(5,2){\vector(1,0){4}}  
\put(7,0){\vector(0,1){4}}  

\put(10,2){\vector(1,0){4}}  
\put(12,0){\vector(0,1){4}}  

\put(2,2){{\color{red}\circle*{0.25}}}
\put(7,2){{\color{red}\circle*{0.25}}}
\put(12,2){{\color{red}\circle*{0.25}}}

\put(6,2){{\color{red}\circle*{0.25}}}
\put(8,2){{\color{red}\circle*{0.25}}}

\put(11.5,2){{\color{red}\circle*{0.25}}}
\put(12.5,2){{\color{red}\circle*{0.25}}}

\linethickness{1mm}
\put(0,2){\color{red}\line(1,0){1}}
\put(3,2){\color{red}\line(1,0){.8}}

\put(5,2){\color{red}\line(1,0){1}}
\put(8,2){\color{red}\line(1,0){.8}}

\put(10,2){\color{red}\line(1,0){1}}
\put(13,2){\color{red}\line(1,0){.8}}

\put(.8,.3){\kern3pt $a <0$}
\put(5.8,.3){\kern3pt $a =0$ (cubic NLS)}
\put(10.8,.3){\kern3pt $a>0$}

\end{picture}

\caption{The spectrum of $\mathcal{H}_\omega$ for the cubic-quintic \eqref{NLS}, namely $F(z) = z^2 + a z^3$ after \cite{PKA}. Spectrum is in red, lines represent essential spectrum and points discrete spectrum or edge resonances. For $a < 0$, there are neither internal modes nor resonances; for $a=0$ (cubic case), a resonance appear; it bifurcates into an internal mode for $a>0$. }
\end{figure}

\subsection{Decay estimates} Decay estimates for the group $e^{it \mathcal{H}_\omega}$ are key to the stability of solitary waves. They can only hold on the essential spectrum; we let $P_e$ denote the associated projector. 

In the absence of embedded eigenvalues, dispersive and Schr\"odinger estimates for the group $e^{it\mathcal{H}_\omega}$ are identical to those for $e^{it\partial_x^2}$, see \cite{BP,KS}
\begin{align*}
& \| P_e e^{it \mathcal{H}_\omega} f \|_{L^{p}} \lesssim |t|^{\frac{1}{p}-\frac{1}{2}} \| f \|_{L^{p'}} \qquad  \mbox{if $2\leq p \leq \infty $} \\
& \| P_e e^{it \mathcal{H}_\omega} f \|_{L^p_t L^q_x} \lesssim \| f \|_{L^2} \qquad  \mbox{if $\frac{1}{p}+\frac{2}{q} = 1$, $p,q>2$}.
\end{align*}

In addition to these classical dispersive and Strichartz estimates, it may come as a surprise that the group $e^{it\mathcal{H}_\omega}$ enjoys stronger local estimates than the flat case $e^{it\partial_x^2}$. Indeed, in the absence of embedded eigenvalues and edge resonances\footnote{Note that $\partial_x^2$ does have an edge resonance at zero frequency, since $\partial_x^2 \, 1 =0$}, we have the improved local decay estimate
$$
\| \langle x \rangle^{-1} e^{it \mathcal{H}_\omega} f \|_{L^\infty} \lesssim |t|^{-\frac 32} \|\langle x \rangle f \|_{L^1}.
$$
This is the sharp estimate \cite{KS}, following \cite{BP} and the sharp scalar case \cite{Schlag}. When a resonance is present, improved local decay can be retrieved after projecting off the resonance \cite{Goldberg2007,Li2024}.

\subsection{Modulation and application to the nonlinear problem} 
Though the decay estimates give a mechanism for stability in the span of the essential spectrum, perturbations in the span of the discrete spectrum will not decay at the linear level. Let us assume for simplicity that no internal modes are present, so that the discrete spectrum reduces to the generalized kernel. We now explain how modulation is the key to dealing with the discrete spectrum.

The idea is to add time-dependent parameters $(\omega,p,\gamma,y)$ which span the manifold of solitary waves. More specifically, we choose the ansatz
$$
v(t, x) = e^{i(p(t)x - \gamma(t))} \left[ (\Phi_{\omega(t)} + u(t))(x + y(t)) \right].
$$

Complexifying the equation as earlier by setting $U = \begin{pmatrix} u \\ \overline{u} \end{pmatrix}$, and plugging the above ansatz, the equation becomes
\begin{equation}
\label{equationperturbation}
i \partial_t U + \mathcal{H}_{\overline{\omega}} U = \underbrace{V U^2 + \dots}_{\mbox{expansion of the nonlinearity} } + \underbrace{ i \dot \omega(t) \Xi_1 + \dots}_{\mbox{modulation terms} } ,
\end{equation}
where $\Xi_1 = \begin{pmatrix} \partial_\omega \Phi \\ \partial_\omega \Phi \end{pmatrix}$, and we omitted most terms on the right-hand side for the sake of clarity. Here arises a subtle point which we will mention and then gloss over. Namely, we linearize around the soliton $\Phi_{\overline{\omega}}$ with parameter $\overline{\omega}$, which is the soliton parameter "at $t=\infty$", and which does not, in general, agree with the current soliton parameter $\omega(t)$...

The coefficients $(\omega,p,\gamma,y)$ are now chosen so that $P_0 U = 0$, where $P_0$ is the spectral projector on the generalized kernel. To obtain an evolution equation for these parameters, we project \eqref{equationperturbation} on the generalized kernel. This ultimately results into an equation of the type
$$
\frac{d}{dt} \left(p(t),\gamma(t),\omega(t),y(t) \right) = \{ \mbox{nonlinear terms} \},
$$
where we need to show that the nonlinear terms decay sufficiently fast for the parameters to converge.

Similarly, projecting \eqref{equationperturbation} on the essential spectrum gives
$$
i \partial_t P_e U + \mathcal{H}_{\overline{\omega}} P_e U = P_e \left[ V U^2 + \{ \mbox{higher-order terms} \} + \{ \mbox{modulation terms} \} \right]
$$

The worst term in the above right-hand side is $V U^2$, and other terms are omitted. It is now easy to understand why the improved local decay estimate is key to the stability of solitons, 
The point is that the term $V U^2$ is not time integrable in the absence of improved local decay: if $U$ decays like $t^{-\frac 12}$ in $L^\infty$, then $UV^2$ has size $\sim t^{-1}$ in $L^2$, which is not time integrable. But improved local decay gives a much stronger control of $V U^2$, and ultimately allows to control the nonlinear problem by resorting only to the various decay estimates stated above (dispersive, Strichartz, improved local). One of the main remaining difficulties is the control of the weighted norm on the profile $f$, which is needed to obtain improved local decay.

\subsection{Full asymptotic stability via decay} By following the approach sketched above, Buslaev and Perelman proved the following theorem, which was the first such result in dimension one.

\begin{theorem}[Buslaev-Perelman \cite{BP} - simplified statement] Assume that $F(z) = O(z^5)$, that $c_\omega>0$ and that the linearized operator does not have embedded eigenvalues, internal modes or edge resonances.
If furthermore $\| v_0 - \Phi_{\omega_0} \|_{H^1 \cap L^{2,2}}$ is sufficiently small, then there exists parameters $(p,\omega,\theta,y)$ and a scattering state $f_\infty$ such that
$$
v(t) = e^{i( p(t)x + \theta(t))} \Phi_{\omega(t)} (x + y(t)) + e^{it\partial_x^2} f_\infty + o_{L^2}(1) \qquad \mbox{as $t \to \infty$}.
$$
\end{theorem}

The above theorem was the first to address asymptotic stability of stationary waves in dimension $1$. Its main drawback is the requirement that $F$ vanishes at order $5$: at the level of \eqref{NLS}, this means that the nonlinearity has order $9$ at least. 

The following result by Krieger and Schlag used sharper local decay estimates and Strichartz estimates, thus weakening the hypothesis to a vanishing of $F$ at order $3+$; furthermore, it provided a precise justification of many steps which were only sketched in \cite{BP}. It focused on the pure power case $F(z) = z^{\sigma+1}$ with $\sigma>2$; as is depicted in Figure \ref{spectrum-p}, unstable eigenvalues are present and full asymptotic stability can only be obtained on a stable manifold. We also mention \cite{Mizumachi} which used similar tools to deal with the simpler problem of small solitary waves generated by an exterior potential; and \cite{MasakiMurphySegata2,MasakiMurphySegata3} where the case of a Dirac potential is considered.

\begin{theorem}[Krieger-Schlag \cite{KS} - simplified statement] \label{theoKS} Assume that $F(z) = z^{\sigma+1}$ with $\sigma >2$. If $v_0$ is sufficiently close to $\Phi_{\omega_0}$ for the norm $\| f \|_{H^1} + \| \langle x \rangle f \|_{L^1 \cap L^2} + \| \langle x \rangle \partial_x f \|_{L^1}$, then full asymptotic stability holds on a manifold of codimension one. Namely, on that manifold, there exists parameters $(p,\omega,\theta,y)$ and a scattering state $f_\infty$ such that
$$
v(t) = e^{i( p(t)x + \theta(t))} \Phi_{\omega(t)} (x + y(t)) + e^{it\partial_x^2} f_\infty + o_{L^2}(1) \qquad \mbox{as $t \to \infty$}.
$$
\end{theorem}

\section{Distorted Fourier transform and nonlinear resonances}

By exploiting only the decay of the linearized problem, it seems impossible to prove full asymptotic stability for potentials $F$ vanishing to order $< 3$, and Theorem \ref{theoKS} seems optimal in this respect. 

But generic and physically relevant potentials vanish to order $2$! In order to deal with such potentials, it is necessary to take advantage of the time oscillations induced by the linearized operator, which play a key role in the nonlinear problem through resonances. This section is dedicated to the exploration of this idea.

\subsection{The distorted Fourier transform} In order to see most clearly the resonances induced by the linearized operator $\mathcal{H}_\omega$, it is convenient to diagonalize this operator, which is the object of the distorted Fourier transform. The linear operator $\mathcal{H}_\omega$ is vectorial, making details significantly more intricate than for scalar Schr\"odinger operators, for which the corresponding theory is classical \cite{DeiftTrubowitz,Yafaev}.

We will now summarize the theory of the distorted Fourier transform, referring to \cite{BP,CollotGermain,KS} for further details. The distorted Fourier transform of a function from $\mathbb{R}$ to $\mathbb{C}^2$ is defined as the pair of scalar functions $(\widetilde{f}_+,\widetilde{f}_-)$ given by (once again omitting indices)
$$
\widetilde{f}_{\pm} (\xi) = \frac{1}{\sqrt{2\pi}} \int f(x) \cdot \sigma_3 \overline{ \psi_{\pm}(x,\xi)} \dd x, \qquad \mbox{for $\xi \in \mathbb{R}$}
$$
where $\sigma_3$ is the Pauli matrix $\begin{pmatrix} 1 & 0 \\ 0 & -1 \end{pmatrix}$ and the functions $\psi_{\pm}(x,\xi)$ are bounded solutions of 
$$
\mathcal{H}_\omega \psi_{\pm}(x,\xi) = \pm (\omega + \xi^2) \mathcal{H}_\omega \psi_{\pm}(x,\xi)
$$
The projection of the function $f$ on the essential spectrum can be recovered through the formula
$$
P_e f(x) = \frac{1}{\sqrt{2\pi}} \sum_\pm \pm \int \widetilde{f}_\pm(\xi) \psi_\pm(x,\xi) \,d\xi
$$
(analogous to the inverse Fourier transform). 
Finally, denoting $\widetilde{\mathcal{F}}$ for the map $f \mapsto (\widetilde{f}_+,\widetilde{f}_-)$, it diagonalizes $\mathcal{H}_\omega$ in the sense that
$$
\widetilde{\mathcal{F}}^{-1} (\pm) (\omega + \xi^2) \widetilde{\mathcal{F}} = \mathcal{H}_\omega.
$$
(the analog for the classical Fourier transform $\widehat{\mathcal{F}}$ being the formula $\widehat{\mathcal{F}}^{-1} \xi^2 \widehat{\mathcal{F}} = -\partial_x^2$).

While the distorted transform $\widetilde{\mathcal{F}}$ has much in common with its classical cousin, there are also important differences. First, it is not unitary, since $\mathcal{H}_\omega$ is not self-adjoint. More importantly for our purposes, the operator $\mathcal{H}_\omega$ is not translation invariant, and the formula $\widehat{fg} = \widehat{f} * \widehat{g}$ is lost. This formula is crucial to treat nonlinear terms after the Fourier transform has been taken; what can replace it will be the subject of the next subsection.

\subsection{Nonlinear spectral distributions} With the help of the distorted Fourier transform, we can diagonalize the linear part of the equation
$$
i \partial_t U + \mathcal{H}_\omega U = 0,
$$
which gives an essentially complete understanding of the linear dynamics... But there remains the nonlinear terms, which we will also need to view on the distorted Fourier side. Amongst the many terms which arise when writing the equation for the perturbation $U$, the most significant are $V(U,U)$ (where $V$ is a smooth and decaying matrix potential) and $\begin{pmatrix} U_1^2 U_2 \\ -  U_2^2 U_1 \end{pmatrix}$, which corresponds to the term $|u|^2 u$ written in $U$ coordinates.

We can express the distorted Fourier transforms of these terms through the nonlinear spectral distributions $\mu$ and $\nu$:
\begin{align*}
& \widetilde{\mathcal{F}}(V(g,g))(\xi) = \int \widetilde{g}(\eta) \widetilde{g}(\sigma) \mu(\xi,\eta,\sigma) \dd \eta \dd \sigma \\
& \widetilde{\mathcal{F}} \begin{pmatrix} g_1^2 g_2 \\ -  g_2^2 g_1 \end{pmatrix} = \int \widetilde{g}(\eta) \widetilde{g}(\sigma) \widetilde{g}(\zeta) \nu(\xi,\eta,\sigma)  \dd \eta \dd \sigma.
\end{align*}

Here, we dropped all indices on the right-hand side to make notations lighter - but $g$, $\widetilde{g}$ take values in $\mathbb{C}^2$ while $\mu$ and $\nu$ take values in $\mathbb{C}^{2 \times 2 \times 2}$ and $\mathbb{C}^{2 \times 2 \times 2 \times 2}$ respectively. Since the matrix potential $V$ is in the Schwartz class, the associated quadratic spectral distribution $\mu$ is smooth and localized. But the cubic term does not feature any localized potential, and the associated cubic spectral distribution $\nu$ is given by
$$
\nu(\xi,\eta,\sigma.\zeta) =\sum \delta(\xi \pm \eta \pm \sigma \pm \zeta)\times \{ \mbox{smooth} \} + \sum \operatorname{p.v.} \frac{1}{\xi \pm \eta \pm \sigma \pm \zeta} \times \{ \mbox{smooth} \} +  \{ \mbox{smooth} \}.
$$
The sums above are over all possible choices of signs, and $\{ \mbox{smooth} \}$ stands for a smooth function.

That this formula can be thought of as a generalization of the classical formula $\widehat{fg} = \frac{1}{\sqrt{2\pi}} \widehat{f} * \widehat{g}$ becomes obvious after writing the convolution of $\widehat{f}$ and $\widehat{g}$ under the form
$$
\widehat{f} * \widehat{g} (\xi) = \int \widehat{f}(\eta) \widehat{f}(\sigma) \delta (\xi - \eta - \sigma) \dd\eta \dd \sigma
$$

\subsection{Nonlinear resonances}
Writing $U$ for the perturbation to the soliton (recall that it is valued in $\mathbb{C}^2$ and stands for $\begin{pmatrix} u \\ \overline{u} \end{pmatrix}$), the equation we are facing is
$$
i \partial_t U + \mathcal{H}_{\overline{\omega}} U = V(U,U) + c \begin{pmatrix} U_1^2 U_2 \\ -  U_2^2 U_1 \end{pmatrix} + \{ \mbox{higher order terms} \} + \{ \mbox{modulation terms} \} + \dots.
$$
Here, $V$ is a decaying matrix potential which can be expressed in terms of $F$ and $\Phi$; as for $\begin{pmatrix} U_1^2 U_2 \\ -  U_2^2 U_1 \end{pmatrix}$, it corresponds to the term $|u|^2 u$ written in $U$ coordinates, and it is present since we are assuming that $F$ is quadratic to leading order.

We want to proceed in the same way that we did when deriving the equation \eqref{equationherisson} and take the distorted Fourier transform before examining resonances. To carry out this plan, we need to understand the terms $V(U,U)$ and $\begin{pmatrix} U_1^2 U_2 \\ -  U_2^2 U_1 \end{pmatrix}$ when viewed in distorted Fourier space - the other terms as well, but these two are the most problematic!

We now change the unknown function to work with to the ($\mathbb{C}^2$-valued) profile
$$
f = e^{-it \mathcal{H}_{\overline{\omega}}} U.
$$
Using the nonlinear spectral distributions introduced in the previous subsection to write the equation on $\widetilde{f}$, it takes the form
\begin{equation}
\label{eqftilde}
\begin{split}
& i \widetilde{f}(\xi) = \int_0^t \int e^{it \Phi_2(\xi,\eta,\sigma)} \widetilde{f}(\eta) \widetilde{f}(\sigma) \mu(\xi,\eta,\sigma) \dd \eta \dd \sigma \dd t + \int_0^t \int e^{it \Phi_3(\xi,\eta,\sigma,\zeta)} \widetilde{f}(\eta) \widetilde{f}(\sigma) \widetilde{f}(\zeta) \nu(\xi,\eta,\sigma)  \dd \eta \dd \sigma \dd\zeta \dd t\\
& \qquad \qquad \qquad \qquad \qquad + \{ \mbox{higher order terms} \};
\end{split}
\end{equation}
as above, this is a simplified expression where indices have been omitted. The quadratic and cubic phases are
\begin{align*}
& \Phi_2(\xi,\eta,\sigma) = [ \xi^2 + \overline{\omega} ] \pm [ \eta^2 + \overline{\omega} ] \pm [ \sigma^2 + \overline{\omega} ] \\
& \Phi_3(\xi,\eta,\sigma,\zeta) = [ \xi^2 + \overline{\omega} ] \pm [ \eta^2 + \overline{\omega} ] \pm [ \sigma^2 + \overline{\omega} ] \pm [ \zeta^2 + \overline{\omega} ],
\end{align*}
 and the signs are $+$ or $-$ depending on the various indices which were omitted.

The idea is now to regard the two integrals on the right-hand side as oscillatory integrals. From this point of view, several features will play an important role
\begin{itemize}
\item Time resonances correspond to points at which the phase $\Phi$ (either $\Phi_2$ or $\Phi_3$) vanishes. Interactions of frequencies for which $\Phi$ is not zero will be weakened through averaging by time oscillations.
\item Space resonances correspond to points at which the gradient of the phase $\nabla \Phi$ (in the frequency variables $\eta,\sigma,\zeta$) vanishes. Away from these points, time decay will follow from oscillations in the frequency variables.
\item Space-time resonances are points at which $\Phi = \nabla \Phi =0$, and they are the main obstructions to decay of $U$, as explained in the two previous points.
\item We also need to take into account the nonlinear spectral distribution, and more specifically its singular set. For quadratic terms, it is empty, and this does not play any role; but for cubic terms, it is of the type $\xi \pm \eta \pm \sigma \pm \zeta = 0$, and oscillations of the phase only matter along this singular set, rather than across.
\item Finally, this whole discussion is only justified as long as the profiles $f$ are sufficiently well-controlled: we have to show as part of the analysis that they vary slowly with time, and that they maintain some regularity in the frequency variable - typically, one propagates an $H^1$ norm in frequency.
\end{itemize}

The method which was sketched, relying on the formalism of oscillatory integrals, is an adaptation of the method of space-time resonances \cite{Germain,GMS1,GMS2}, to the context of asymptotic stability of solitary waves. An intermediary step consisted of the study of cubic \eqref{NLS} with exterior potential, which was pursued by a number of authors \cite{CP1,CP2,Delort,GPR,MasakiMurphySegata,Naumkin,Naumkin2,Stewart}.

We focused in this subsection on the control of the radiation, but modulation terms must also be taken into account. For them too, nonlinear resonances must be tracked since decay of the radiation is not enough to bound them.

\subsection{Full asymptotic stability via nonlinear resonances}
The following theorem can be proved by the approach that was sketched above; it was preceded by \cite{Chen}, which deals by the same means with the case of small solitary waves generated by exterior potentials.

\begin{theorem}[Collot-Germain \cite{CollotGermain} - simplified statement] Assume that $c_\omega>0$ and that the linearized operator does not have embedded eigenvalues, internal modes or edge resonances.
If furthermore $\| v_0 - \Phi_{\omega} \|_{H^1 \cap L^{2,1}}$ is sufficiently small, then there exists parameters $(p,\omega,\theta,y)$ and a solution $u_{MS}$ undergoing modified scattering such that
$$
v(t) = e^{i( p(t)x + \theta(t))} \Phi_{\omega(t)} (x + y(t)) + u_{MS} + o_{L^2}(1) \qquad \mbox{as $t \to \infty$}.
$$
\end{theorem}

This theorem gives a satisfactory answer to the question of full asymptotic stability if neither edge resonances nor internal modes are present. When they are present, the problem is much more challenging; this will be the subject of the next two sections.

\section{Edge resonances and the cubic (NLS)} 

\label{sectionresonances}

\subsection{The problem with edge resonances} \subsubsection{Quadratic interactions} To understand the specifif problem posed by edge resonances, we consider first the quadratic term appearing in \eqref{eqftilde}
\begin{equation}
    \label{quadratic}
\int_0^t \int  e^{it \Phi_2(\xi,\eta,\sigma)} \widetilde{f}(\eta) \widetilde{f}(\sigma) \mu(\xi,\eta,\sigma) \dd \eta \dd \sigma \dd t 
\end{equation}
with 
$$
 \Phi_2(\xi,\eta,\sigma) = [ \xi^2 + \overline{\omega} ] - [ \eta^2 + \overline{\omega} ] - [\sigma^2 + \overline{\omega} ] = \xi^2 - \eta^2 - \sigma^2 - \overline{\omega}
$$
(notice the specific choice of signs in $\Phi_2$). We now compute space-time resonances, in other words frequencies $(\xi,\eta,\sigma)$ such that
$$
\partial_\eta \Phi_2(\xi,\eta,\sigma) = \partial_\sigma \Phi_2(\xi,\eta,\sigma) = \Phi_2(\xi,\eta,\sigma) = 0 \quad \Longleftrightarrow \quad (\xi,\eta,\sigma) = (\pm \overline{\omega}^{\frac 12} , 0,0).
$$

This specific interaction is degenerate in the sense that the phase is stationary in all the integration variables, namely $\eta,\sigma,t$... But the nonlinear spectral distribution might help! Indeed, in the absence of an edge resonance, we have $\mu(\xi,0,0) = 0$ for any $\xi$ (this is essentially equivalent to improved local decay), which effectively cancels the space-time resonance.

If there is an edge resonance however, such a strong cancellation is not present. As a result of the space-time resonance which we pointed out, logarithmic growth will occur for the profile $\widetilde{f}$ at the frequency $\xi = \pm \overline{\omega}^{\frac 12}$, as observed by Lindblad-L\"uhrman-Schlag-Soffer \cite{LLSS}. Translating this in terms of the solution $U = e^{it \mathcal{H}_\omega} f$ itself, these same authors find that
$$
\| U(t) \|_{\infty} \sim \frac{\log t }{t^{\frac 12}} \qquad \mbox{as $t \to \infty$}.
$$

This is only a logarithmic loss  compared to the dispersive linear decay, but this will cause our method of proof to break down.

\subsubsection{Cubic interactions} The difficulty, which has not been overcome for the time being, stems from cubic interactions involving the "bad" frequencies $\pm \overline{\omega}^{\frac{1}{2}}$.
Arguing heuristically first, we observe that the logarithmic loss in decay  pointed out above has important consequences when we consider the cubic term. Indeed, it was already critical as noted in \eqref{loggrowth}, and any slower decay will have sizeable effects.

To be more precise, we now consider the cubic term 
$$
 \int_0^t \int e^{it \Phi_3(\xi,\eta,\sigma,\zeta)} \widetilde{f}(\eta) \widetilde{f}(\sigma) \widetilde{f}(\zeta) \nu(\xi,\eta,\sigma)  \dd \eta \dd \sigma \dd\zeta \dd t
$$
with the following choice of signs in the cubic phase
$$
 \Phi_3(\xi,\eta,\sigma,\zeta) = [ \xi^2 + \overline{\omega} ] - [ \eta^2 + \overline{\omega} ] + [ \sigma^2 + \overline{\omega} ] - [ \zeta^2 + \overline{\omega} ] =  \xi^2 -  \eta^2 + \sigma^2 - \zeta^2.
$$
For the right choice of signs, the singular set of the nonlinear spectral distribution is $\{ \xi  -  \eta + \sigma - \zeta = 0 \}$. With this singular set and the phase above, we are essentially brought back to \eqref{equationherisson}. In particular, the frequency interaction $(\xi,\eta,\sigma,\zeta) = (\overline{\omega}^{\frac 12}, \overline{\omega}^{\frac 12}, \overline{\omega}^{\frac 12}, \overline{\omega}^{\frac 12})$ is a space-time resonant. This self-interaction of the bad frequecy $\omega^{\frac 12}$ is very hard to control at the nonlinear level.

In the essentially equivalent context of nonlinear Klein-Gordon, L\"uhrmann and Schlag \cite{LuhrmannSchlag2} were able to push the argument up to exponential time scales in the inverse of the size of the perturbation. Going beyond this threshold and obtaining a global result seems very hard! One way to explain the difficulty is that the correct ansatz for $\widetilde{f}$ is unclear: it should still involve modified scattering (phase correction), and account for the logarithmic slowdown... but the precise form of the correction is difficult to pin down.

\subsection{Cubic (NLS) and integrability}
The cubic \eqref{NLS} 
$$
i \partial_t v - \partial_x^2 v - |v|^2 v = 0
$$
is probably the most important example in the class of nonlinear Schr\"odinger equations. As noted in \cite{CGNT}, the linearized operator for the soliton features an edge resonance. As a consequence, the methods discussed above seem unable to prove full asymptotic stability.

However, it was discovered by Zakharov and Shabat \cite{ZakharovShabat} that the cubic \eqref{NLS} is completely integrable, and explicitly solvable through the Inverse Scattering Transform. This approach provides a means of showing full asymptotic stability of the solitary wave, which was exploited by Cuccagna-Pelinovsky \cite{CP}, and extended to multi-soliton configurations  \cite{BJM,Saalmann}.

\subsection{The null form in cubic (NLS)} Surprisingly, it is possible to prove full asymptotic stability for the cubic \eqref{NLS} equation without resorting to the completely integrable structure, in spite of the seemingly unsurmountable difficulty explained above.

The key is to look more precisely into the structure of the quadratic term \eqref{quadratic}. We saw that the worst interaction was the space-time resonance at frequencies $(\overline{\omega}^{\frac 12},0,0)$; and we also saw that the nonlinear spectral distribution $\mu(\xi,\eta,\sigma)$ vanishes for $\eta = \sigma = 0$ in the absence of an edge resonance. The crucial observation of Li and L\"uhrmann \cite{LiLuhrmann24} is that, in the case of cubic \eqref{NLS}, $\mu(\overline{\omega}^{\frac 12},0,0) =0$, as follows from an explicit computation! While this observation allows to deal with the radiation part of the solution, it is not sufficient to obtain control of the modulation; and indeed a second cancellation can be identified in the modulation terms! Leveraging these two cancellation allowed these authors to prove the following theorem.

\begin{theorem}[Li-L\"uhrmann \cite{LiLuhrmann24} - simplified statement] Consider the cubic \eqref{NLS} with even data. If $\| v_0 - \Phi_{1} \|_{H^1 \cap L^{2,1}}$ is sufficiently small, then there exists parameters $(\omega,\theta)$ and a solution $u_{MS}$ undergoing modified scattering such that
$$
v(t) = e^{i \theta(t)} \Phi_{\omega(t)} (x) + u_{MS} + o_{L^2}(1) \qquad \mbox{as $t \to \infty$}.
$$
\end{theorem}

\section{Internal modes and radiation damping} 

\label{sectioninternal}

\subsection{The mechanism of radiation damping} In the presence of an internal mode, the linearized equation would not lead to asymptotic stability - indeed, the internal mode would oscillate indefinitely without decaying. 

A striking realization is that the coupling with the continuous spectrum (under appropriate assumptions) leads to the decay of the internal mode. To be more specific, the nonlinearity enables energy to be transfered from the discrete to the continuous spectrum, and then carried away to infinity through dispersive decay. The first mathematical proof of this phenomenon is due to Soffer and Weinstein \cite{SofferWeinstein}.

To explain very roughly the mechanism of radiation damping, we will henceforth ignore modulation questions and decompose the perturbation of the solitary wave into
$$
U(t) = a(t) \begin{pmatrix} \varphi \\ \varphi \end{pmatrix} + U_e (t) + \{ \mbox{generalized kernel elements} \},
$$
where $\begin{pmatrix} \varphi \\ \varphi \end{pmatrix}$ is the internal mode and $U_e$ is the projection of $U$ onto the essential spectrum of the operator.

The evolution equation for $U$ can then be written as a coupled system for $a$ and $U_e$:
\begin{align*}
& i \partial_t a - i \lambda a = P_i (V \overline{a} \varphi U_e) + \dots \\
& i \partial_t U_e + \mathcal{H}_\omega U_e = P_e (V a^2 \varphi^2) + \dots .
\end{align*}
The notations used here are quite loose; in particular, the vectorial nature of the problem is mostly ignored. Most terms are omitted on the right-hand side: the only terms that are kept are the ones responsible for radiation damping. Finally, $P_i$ and $P_e$ stand for the projections on the internal mode and the essential spectrum respectively, and $\lambda$ is the time frequency associated to the internal mode.

From the ODE satisfied by $a$, we deduce after setting $A(t) = a(t) e^{i\lambda t}$ that
$$
i \partial_t A(t) = e^{2i \lambda t} P_i ( V \overline{A} \varphi U_e) + \dots 
$$
Turning to the equation satisfied by $U_e$, we substitute $A$ for $a$,  write the Duhamel formula and use a rough approximation to obtain
\begin{align*}
i U_e(t) & = \int_0^t e^{i(t-s) \mathcal{H}_\omega} A(s)^2 e^{-2i\lambda s} \dd s \, P_e (V \varphi^2) \\
& \sim A(t)^2 \int_0^t e^{i(t-s) \mathcal{H}_\omega} e^{-2i\lambda s} \dd s \, P_e (V \varphi^2) \\
& \sim e^{-2i\lambda t} A(t)^2 \delta(\mathcal{H} + 2 \lambda) P_e (V \varphi^2).
\end{align*}

Combining the equations on $\partial_t A$ and $U_e$ gives
$$
\partial_t A(t) = - |A|^2 A \, P_i(V \varphi \delta(\mathcal{H} + 2 \lambda) P_e (V \varphi^2)).
$$
Provided the coefficient $P_i(V \varphi \delta(\mathcal{H} + 2 \lambda) P_e (V \varphi^2))$ is positive, this gives the desired decay of $A(t)$. This condition is known as \textit{nonlinear Fermi golden rule}.

\subsection{Full asymptotic stability through radiation damping} By exploiting the radiation damping mechanism explained above, Buslaev and Sulem were able to show the following theorem.

\begin{theorem}[Buslaev-Sulem \cite{BS} - simplified statement] Assume that $F(z) = O(z^5)$, that $c_\omega >0$, and that the linearized operator has one internal mode, but no embedded spectrum or edge resonances.

If furthermore the data $v_0$ is even and $\| v_0 - \Phi_{\omega} \|_{H^1 \cap L^{2,2}}$ is sufficiently small, then there exists parameters $(p,\omega,\theta,y)$ and a scattering state $f_\infty$ such that
$$
v(t) = e^{i( p(t)x + \theta(t))} \Phi_{\omega(t)} (x + y(t)) + \{\mbox{time-decaying solution}\} \qquad \mbox{as $t \to \infty$}.
$$
\end{theorem}

Though this theorem shows that nonlinear stability results from radiation damping, it does so under a very strong assumption on the nonlinearity, namely that it vanishes to order 9 or higher... It would of course be desirable to relax this condition.

The problem one faces to do so is similar to the one occuring in the presence of edge resonances. Namely, the radiation damping mechanism results in the growth of a given (distorted) frequency, and the self-interactions of this frequency through the cubic term (in the case where the nonlinearity is cubic to leading order, say) seem very difficult to control. In the related context of nonlinear Klein-Gordon equations, the best result in this direction is due to Delort-Masmoudi \cite{DelortMasmoudi}, who were able to reach a time scale $\sim \epsilon^{-4}$, denoting $\epsilon$ the size of the initial perturbation.

Given the paucity of results in dimension 1, it is worthwhile looking into the case of dimension 3, where the picture is much more complete. The seminal work of Soffer-Weinstein \cite{SofferWeinstein} on radiation damping for cubic Klein-Gordon equations was extended by \cite{BambusiCuccagna,TsaiYau}. The case of quadratic Klein-Gordon seemed out of reach until the very recent work by L\'eger and Pusateri \cite{LegerPusateri,LegerPusateri2022} who analyzed nonlinear resonances with the help of harmonic analytic tools developed in \cite{PusateriSoffer}.

\section{The virial method and local asymptotic stability}

We have been reviewing up to this point the spectral approach to full asymptotic stability; in this section, we turn to the virial approach. It was first applied to nonlinear Klein-Gordon equations (see the next section) and has proved extremely successful. This approach relies on a very different point of view on the question of asymptotic stability, with complementary results under complementary assumptions to the spectral approach.
To give a meaningful introduction to the virial method in the framework of \eqref{NLS} would go beyond the scope of this review, and we refer to \cite{Martel2024}. Rather, we will try here to compare the virial and spectral approaches as far as their aims, their tools, their outcomes and their shortcomings go.

\begin{itemize}
\item \textit{Local versus full asymptotic stability} (these two notions are defined in the introduction). The former can be proved through the virial method and the latter through the spectral method.

\smallskip 

\item \textit{Fourier versus energy methods.} As is already clear from their names, the heart of the virial method consists of energy estimates while the spectral method relies on estimates in Fourier or spectral space. As such, the virial method can tolerate more nonlinear situations while the spectral method is more closely tied to linear(ized) behavior.  

\smallskip 

\item \textit{Energy space versus weighted spaces.} Since it relies on energy estimates, the virial method naturally applies to perturbations in the energy space $H^1$. As far as the spectral method goes, it seems to require perturbations in weighted spaces; this condition is in particular crucial to obtain local improved decay and to take advantage of nonlinear resonances.

\smallskip 

\item \textit{Small solitary waves.} By its very definition, full asymptotic stability excludes additional small solitary waves in the limit $t \to \infty$. By contrast, they are not precluded by the virial method; indeed, taking the perturbation to be small in the energy space does not exclude them, see Subsection \ref{subsectionsolitary}. Since these small solitary waves are asymptotically moving away from the main solitary wave, they are not violating the definition of local asymptotic stability.

\smallskip 

\item \textit{Breathers.} For the cubic \eqref{NLS} and other integrable models, breathers exist in an $H^1$ neighborhood of solitary waves, see Subsection \ref{subsectionbreathers}; thus, they need to be taken into account when applying the virial method, which has proved difficult. Weighted topologies exclude breathers, which for instance do not appear in the theorem of Li and L\"uhrmann stated earlier.

\smallskip 

\item \textit{Decay rate.} In the known applications of the spectral method, it is essential to be able to control the decay of the radiation in an appropriate norm by an inverse power of $t$. This is not the case for the virial method, and indeed, the presence of small solitary waves shows that a decay rate for the radiation in a window around the solitary wave cannot be expected.

\smallskip 

\item \textit{Limitations.} For reasons which have been explained above, full asymptotic stability in the presence of edge resonances or internal modes seems currently out of reach of the spectral method. As for the virial method, it was applied successfully to deal with internal modes, but edge resonances seem to constitute an obstruction.
\end{itemize}

The application of the virial method to asymptotic stability in \eqref{NLS} is very recent. The theorem by Martel stated below was seminal, and it immediately sparked further developments \cite{CuccagnaMaeda1, CuccagnaMaeda2,Martel1, Martel2, Rialland, Rialland2}. 

\begin{theorem}[Martel \cite{Martel1,Martel2}] In the cubic-quintic case $F(z) = z^2 + a z^3$, with $a \not 0$, consider a solitary wave $\Phi_\omega$ with $\omega$ sufficiently small. If $v_0$ is even and $\| v_0 - \Phi_\omega \|_{H^1}$ is sufficiently small, then there exists parameters $(\omega,\theta)$ such that
$$
\sup_{N>0} \lim_{t \to \infty} \| v(t) - e^{i \theta(t)} \Phi_{\omega(t)} \|_{L^\infty(-N,N)} = 0.
$$
\end{theorem}

\section{The nonlinear Klein-Gordon equation}

This section will be dedicated to a quick overview of full asymptotic stability for the nonlinear Klein-Gordon equation. We will emphasize how much the development of the theory was parallel to that for \eqref{NLS}: indeed, the questions are remarkably close, and the same tools apply. There are nevertheless differences, which will be highlighted.

\subsection{The equation}
Consider the one-dimensional nonlinear Klein-Gordon equation
\begin{equation}
\label{NLKG} \tag{NLKG}
\begin{cases}
& \partial_t^2 w - \partial_x^2 w - G'(w) = 0 \\
& w(t=0) = w_0 \\
& \partial_t w(t=0) = w_1,
\end{cases}
\end{equation}
where $w(t,x)$ is real-valued, and $(t,x) \in \mathbb{R} \times \mathbb{R}$.
It derives from the energy
$$
\int_{-\infty}^\infty \left[ \frac{1}{2} |\partial_t w |^2 + \frac 12 |\partial_x w|^2 - G(w) \right] \dd x.
$$
We assume that $G''$ does not vanish at the minima of $G$ (to ensure that the resulting PDE, namely \eqref{NLKG}, has a non-trivial mass term).

A stationary solution $\Phi$ is given by the ODE
$$
\partial_x^2 \Phi + G'(\Phi) = 0.
$$
This ODE obviously reduces to the equation \eqref{eq:soliton} for solitary waves of \eqref{NLS} for an appropriate $F$; the condition $\omega>0$ corresponds to $G$ having a non-degenerate local maximum at its minimum. 

The above ODE gives a stationary solution of \eqref{NLKG}, which generates a one-parameter family by space translation. Traveling waves can be deduced by applying the Lorentz transform.

\subsection{Decaying solitary waves} By this, we mean solitary waves $\Phi$ which are localized in space, similar to the solitary waves of \eqref{NLS} that we have been discussing; this corresponds to potentials $G$ with a single well. These solitary waves are never stable, so that asymptotic stability can only be reached on a stable manifold; a more ambitious aim is to classify solutions in a neighborhood of the solitary wave, including solitary waves away from the stable manifold \cite{KriegerNakanishiSchlag}.

Full asymptotic stability was first proved by Krieger-Nakanishi-Schlag \cite{KriegerNakanishiSchlag} via decay estimates, for nonlinearities of order $>5$ (same numerology as \eqref{NLS}). Kairzhan and Pusateri \cite{KairzhanPusateri} analyzed the equation through the distorted Fourier and were able to treat the quartic case.

Local asymptotic stability was obtained in many cases
\cite{CuccagnaMaeda23bis,CuccagnaMaedaScrobogna,KowalczykMartelMunoz2,LiLuhrmann23,PalaciosPusateri}, the main obstacle to the application of this method being the presence of an edge resonance.

\subsection{Topological solitary waves} Solitary waves which have a non-trivial behavior at infinity (nonzero limit) are called topological solitary waves, or kinks in the context of \eqref{NLKG}. Such solitary waves occur if the potential $G$ has a double-well structure, and kinks are connecting the two nimima, converging to either one as $x\to \pm \infty$.

The most classical example is the kink solution $K(x) = \tanh \left( \frac {x}{\sqrt 2} \right)$ of the $\Phi^4$ model given by the potential $G(w) = - \frac{1}{2}(1-w^2)^2$. This was the first model for which local asymptotic stability was established, in the breakthrough paper by Kowalczyk-Martel-Munoz \cite{KowalczykMartelMunoz1}, which was followed by \cite{KowalczykMartelMunoz2017,KMMV,CuccagnaMaeda23,KowalczykMartel}. The virial method introduced by Kowalczyk-Martel-Munoz is very effective to deal with internal modes, but edge resonances remain an obstacle.

Full asymptotic stability was first approached via decay estimates \cite{KK2,KK1} for sufficiently high order nonlinearities. As a next step, toy models were investigated \cite{LLSS,LLS2,LLS1} which incorporated many of the difficulties of the equation for the perturbation of the solitary wave of \eqref{NLKG}. Finally, full asymptotic stability of the kink was obtained in \cite{GP,GPZ} under some restrictions: odd solutions, absence of edge resonances and internal modes. If these spectral assumptions are relaxed, the best known results give stability over large time intervals \cite{DelortMasmoudi,LuhrmannSchlag2}.

Finally, the Sine-Gordon equation, which corresponds to the potential $G(w) = - \cos w$, is known to be integrable. There is a striking analogy with the cubic \eqref{NLS} equation: not only are both models integrable, but they both support a variety of (unstable) breather-like solutions which do not fit into the standard asymptotic decomposition as a sum of radiation and a solitary wave. Full asymptotic stability of the kink could be established by integrable means \cite{ChenLiuLu,KochYu} but also by more direct spectral PDE methods \cite{LuhrmannSchlag1}. As was the case for cubic \eqref{NLS}, the latter approach is successful thanks to a subtle cancellation in quadratic interactions.

\section{Perspectives} The two main gaps in our understanding of full asymptotic stability of solitary waves in 1D nonlinear dispersive equations are the cases where resonances or internal modes are present in the linearized operator; this was discussed in Sections \ref{sectionresonances} and \ref{sectioninternal} respectively. It would be most desirable to close these gaps! 

\medskip

A variety of other semilinear dispersive problems are of interest in dimension one, and it is our belief that the methods developed for \eqref{NLS} will apply there.

\medskip

Going beyond this class of equations, the asymptotic stability of solitary waves is very poorly understood, in particular for
\begin{itemize}
\item Discrete problems (set on the lattice)
\item Quasilinear problems (in particular, involving waves in fluids)
\item Higher-dimensional problems (the case of dimension 1 is harder in that the decay is slower, but in dimension $\geq 2$ the geometry of resonances as well as spectral theory become more challenging).
\end{itemize}
It is our hope that some of the methods reviewed here can find applications to solve these problems!

\medskip

Finally, as mentioned earlier, the soliton resolution conjecture would be the ultimate aim in understanding large-time behavior of nonlinear dispersive problems set on Euclidean space...

\bibliographystyle{abbrv}
\bibliography{references}

\subsection*{Acknowledgements} The author was supported by a Wolfson fellowship from the Royal Society and the Simons collaboration on Wave Turbulence. He is grateful to Charles Collot, Jonas L\"uhrmann and Yvan Martel for their comments on an earlier version of this article.
\end{document}